\documentclass[a4paper,12pt]{amsart}

\usepackage{accents}
\usepackage{bm}
\usepackage{mathtools}

\usepackage{cleveref}

\DeclareMathOperator{\Cov}{Cov}
\DeclareMathOperator{\Dom}{dom}
\DeclareMathOperator{\Expectation}{\mathbb E}
\DeclareMathOperator{\Maxexp}{\mathcal E}
\DeclareMathOperator{\Prob}{\mathbb P}

\newcommand{\Cinfcs}[1]{C_0^\infty\left(#1\right)}
\newcommand{\Cinfp}[1]{C_{\mathrm{p}}^\infty\left(#1\right)}
\newcommand{\KH}[2]{\operatorname{DH}\left(#1 \middle\vert #2 \right)}
\newcommand{\Lexp}[1]{\orliczof {\coshtwo} #1}
\newcommand{\LlogL}[1]{\orliczof {\coshtwo_*} #1}
\newcommand{\M}{\mathcal M}

\newcommand{\Wexptwo}[1]{W^{1,2}_{\coshtwo}\left(#1\right)}
\newcommand{\Wexp}[1]{W^{1,(\cosh-1)}\left(#1\right)}
\newcommand{\WlogL}[1]{W^{1,(\cosh-1)_*}\left(#1\right)}
\newcommand{\Wlogtwo}[1]{W^{1,2}_{\coshtwo_*}\left(#1\right)}
\newcommand{\acceleration}[1]{\accentset{**}{#1}}
\newcommand{\avalof}[1]{\left\vert#1\right\vert}
\newcommand{\coshtwo}{{\cosh_2}}
\newcommand{\covat}[3]{\Cov_{#1}\left(#2,#3\right)}

\newcommand{\cpoly}[2]{C^{#1}_{\text{poly}}\left(\reals^{#2}\right)}
\newcommand{\derivby}[1]{\frac{d}{d#1}}
\newcommand{\displacement}{\operatorname{\mathbb S}}
\newcommand{\domof}[1]{\Dom\left(#1\right)}
\newcommand{\etransport}[2]{\prescript{\text{e}}{} {\mathbb U} _ {#1} ^ {#2}}
\newcommand{\euler}{\mathrm{e}}

\newcommand{\expectat}[2]{\Expectation_{#1}\left[#2\right]}

\newcommand{\expfiberat}[2]{S_{#1}\maxexpat{#2}}
\newcommand{\expof}[1]{\exp\left(#1\right)}
\newcommand{\exptwo}{{\exp_2}}
\newcommand{\gaussdensity}{\gamma}
\newcommand{\gaussint}[2]{\int{#1}\,\gaussdensity(#2)\,d#2 \,}
\newcommand{\gausstwo}{{\operatorname{gauss}_2}}
\newcommand{\integrald}[3]{\int_{#1} {#2} \ {#3}}
\newcommand{\maxexpat}[1]{\Maxexp\left(#1\right)}
\newcommand{\maxexp}[1]{{\mathcal E}\left(#1\right)}
\newcommand{\mixbundleat}[1]{\prescript{*}{}S\maxexpat{#1}}
\newcommand{\mtransport}[2]{\prescript{\text{m}}{} {\mathbb U} _ {#1} ^ {#2}}
\newcommand{\normat}[2]{\left\Vert#2\right\Vert_{#1}}
\newcommand{\normof}[1]{\left\Vert#1\right\Vert}
\newcommand{\orliczof}[2]{L_{#1}\left(#2\right)}
\newcommand{\orliczpof}[3]{L_{#2}^{#1}\left(#3\right)}
\newcommand{\pderivby}[1]{\frac {\partial} {\partial #1}}
\newcommand{\probat}[2]{\Prob_{#1}\left(#2\right)}
\newcommand{\reals}{\mathbb{R}}
\newcommand{\scalarat}[3]{\left\langle#2,#3\right\rangle_{#1}}
\newcommand{\scalarof}[2]{\left\langle#1,#2\right\rangle}

\newcommand{\setof}[2]{\left\{#1 \, \middle\vert \, #2 \right\}}

\newcommand{\transport}[2]{{\mathbb U} _ {#1} ^ {#2}}
\newcommand{\velocity}[1]{\accentset{\star}{#1}}

\title[Statistical bundle on Gaussian Orlicz-Sobolev space]{Affine statistical bundle modeled on a Gaussian Orlicz-Sobolev space}

\author[G. Pistone]{Giovanni Pistone}
\email{giovanni.pistone@carloalberto.org}
\address{de Castro Statistics, Collegio Carlo Alberto, Piazza Vincenzo Arbarello 8, Torino 10122 Italy}
\urladdr{giannidiorestino.it}

\begin{document}

\begin{abstract}The dually flat structure of statistical manifolds can be derived in a non-parametric way from a particular case of affine space defined on a qualified set of probability measures. The statistically natural displacement mapping of the affine space depends on the notion of Fisher's score. The model space must be carefully defined if the state space is not finite. Among various options, we discuss how to use Orlicz-Sobolev spaces with Gaussian weight. Such a fully non-parametric set-up provides tools to discuss intrinsically infinite-dimensional evolution problems. \textbf{keywords} Information geometry, Gaussian Orlicz-Sobolev space, Statistical Bundle, Exponential Manifold, Dually Flat Affine Manifold\end{abstract}

\maketitle

\tableofcontents

\section{Non-parametric statistical bundle}
\label{sec:bundle}
Professor S.-I. Amari clearly stated in a 1987 conference paper \cite{amari:87dual} the notion of a non-parametric fiber bundle in Information Geometry (IG). He says
\begin{quote}
\emph{A fibre bundle is constructed on a finite-dimensional parametric statistical model with a Hilbert space as the fibre space. The Hilbert space represents the tangent directions of the set of probability distributions in the function space. A pair of dual linear connections is introduced in the Hilbert bundle.}
\end{quote}
A related journal paper \cite{amari-kumon:1988-AS} shows the applications to the statistics of semi-parametric statistical models. R. E. Kass and P. W. Vos, in their 1997 monograph \cite[10.3]{kass-vos:1997}, further explain the construction by describing the tangent fiber $T_p\mathcal M$ of a statistical model $\mathcal M$ as a vector space of random variables:
\begin{quote}
 \emph{The tangent space $T_p\mathcal M$ is an inner product space of random variables having zero expectation and finite variance that describe $\mathcal M$ at $p$. Finite-dimensional vector bundles enlarge $T_p\mathcal M$ by adding other random variables, each of which has mean zero and finite variance. The Hilbert bundle carries this process to its conclusion by enlarging $T_p\mathcal M$ to the space of all random variables having zero mean and finite variance.}   
\end{quote}

I started to work on a non-parametric version of professor Amari's ideas in the early nineties \cite{pistone-sempi:95}. In doing that, the self-imposed prescription to avoid parameters led to the use of differential geometry as presented by N. Bourbaki \cite{bourbaki:71variete} and S. Lang \cite{lang:1995}. Specifically, the approach involves the expression of the tangent spaces as vector spaces of random variables. I give below an updated summary of such a construction.  

If $\mu$ is any probability measure, $\mu \in \mathcal M_1$, we define as fiber at $\mu$ the Hilbert space $H_{\mu}$ of $\mu$-square integrable random variable with zero $\mu$-mean, $H_{\mu} = L^2_0(\mu)$. The set of all couples $(\mu,u)$ with $\mu \in \mathcal M_1$ and $u \in H_{\mu}$ is the (maximal) (statistical) Hilbert bundle $H\mathcal M_1$. The restriction to $\mu$-centered random variables is intended to reflect the affine constraint $\mu(X) = 1$ satisfied by $\mathcal M_1$ as a subset of the vector space $\mathcal S$ of finite signed measures.

In which sense does each fiber $H_{\mu}$ contain the tangents to statistical models? Clearly, $M_1$ is a convex set whose tangent space is the set $S_0$ of signed finite measure with total charge 0. But there is another sense to consider. Assume the 1-dimensional statistical model model $\theta \mapsto \mu(\theta)$ is set-wise continuously differentiable, that is, for each measurable set $A \in \mathcal X$ the real map $\theta \mapsto \mu(A;\theta)$ is continuously differentiable with derivative $\theta \mapsto \dot \mu(A;\theta)$. Assume that $\theta \mapsto \dot \mu(\theta)$ is a curve in the Banach space $\mathcal S_0$ of signed finite measures with 0 total charge and total variation norm. Clearly, $\dot\mu(A;\theta) = 0$ if $\dot \mu(A;\theta) = 0$, hence $\dot\mu(\theta)$ is absolutely continuous with respect to $\mu(\theta)$. See this argument in \cite{ay-jost-le-schwachhofer:2017ig} and, in full detail, in \cite{le:2022arXiv}.

The Radon-Nikodym derivative $\ell(\theta) = d\dot\mu(\theta)/d\mu(\theta)$ is such that $\int \ell(\theta) \, d\mu(\theta) = \int d \dot \mu(\theta) = 0$. If moreover $\int \ell(\theta)^2 \, d\mu(\theta) = \int l(\theta) d\dot\mu(\theta)$ is finite, then $\ell(\theta) \in H_\theta = H_{\mu_\theta}$. Under the classical statistical assumption of a regular likelihood, $\mu_\theta = p_\theta \cdot \mu$, $p_\theta > 0$, it holds  \begin{equation*}
\ell(\theta) = \frac{\dot p_\theta}{p_\theta} = \derivby \theta \log p_\theta
\end{equation*}
That is, $\ell_\theta$ is the Fisher's score of the statistical model, a statistically natural expression of the rate of variation \cite[4.2]{efron-hastie:2016}. In this sense, it is an expression of the tangent vector.

Conversely, let be given an element $u \in H_{\mu}$ such that $\int \euler^{\theta u} \, d\mu < \infty$ for all $\theta$ in a neighborhood of 0. The statistical model $\mu(\theta) \propto \euler^{\theta u} \cdot \mu$ is such that the Fisher score at 0 is $u$. This, again, is quite a classical argument. It depends only on the classical theory of exponential families, see \cite{brown:86} or \cite[5.5]{efron-hastie:2016}.

Let us fix a reference measure $\mu$ and restrict our attention to the set $\mathcal P \cdot \mu$, $\mathcal P$ strictly positive probability densities. If the sample space is finite, all conditions outlined above are met. The Hilbert bundle is an expression of the tangent bundle where Fisher's score expresses the velocity of one-dimensional curves. 

If the sample space is not finite, one must introduce restrictions on both the set of probabilities and the random variables in the fiber. One option is to restrict to a convenient subset of positive densities $\mathcal E \subset \mathcal P$ in such a way there is a vector space $B(\mathcal E) \subset \cap_{\mu \in \mathcal E} L^2(\mu)$ that contains the scores of all one-dimensional statistical models. Then, each fiber is defined to be
\begin{equation*}
    S_p \mathcal E = \setof{u \in B(\mathcal E)}{ \int u \, p \, d\mu = 0} \ .
\end{equation*}

For \emph{example}, assume the sample space is a compact set $K$, and $\mu$ is a diffuse measure. Let $\mathcal E$ be the set of all continuous, strictly positive probability densities. $\mathcal E$ is an open convex set of $C(K)$. If $\theta \mapsto p_\theta$ is a differentiable curve in $\mathcal E$, then the score $\derivby \theta \log p_\theta$ is a curve in $B(\mathcal E) = C(K)$ such that for each $\theta$ it holds $\int \derivby \theta \log p_\theta \, d\mu = 0$. And conversely, for each $u \in B_\theta$ the curve $\theta \mapsto p_\theta \propto \euler^{tu}$ has values in $\mathcal E$ and its score at $\theta$ equals $u$. See \cite{bauer-bruveris-michor:2016} for an approach to non-parametric IG based on smooth densities.

A similar but different approach will be used below. Both depend on the idea of finding a convenient class of random variables $u$ such that the exponential family proportional to $\euler^{\theta u}$ is conveniently defined. 

For my purpose, the most important of professor Amari's contributions has been the \emph{definition} a couple of affine connections in a fully nonparametric way, \cite[Th. 1 and 2]{amari:87dual}. In the notations used here, he considers two types of transport between fibers,
\begin{gather}
    \etransport \mu \nu \colon S_\mu \mathcal E \ni u \mapsto u - \int u \, d\nu \in S_\nu \mathcal E \label{eq:etransport} \\
        \mtransport \mu \nu \colon S_\mu \mathcal E \ni  u \mapsto \frac{d\mu}{d\nu} u \, \in S_\nu \mathcal E \label{eq:mtransport}\ ,
\end{gather}
and proves the duality result
\begin{equation*}
   \scalarat \nu {\etransport \mu \nu u} v = \scalarat \mu u {\mtransport \nu \mu v} \ .
\end{equation*}
Also, there is a transport of the inner product from one fiber to the other, 
\begin{equation*}
    \scalarat \mu u v = \scalarat \mu {\etransport \nu \mu \etransport \mu \nu u} v = \scalarat \nu {\etransport \mu \nu u}{\mtransport \mu \nu v} \ .
\end{equation*}

The equations above clearly define a geometry of probability measures that is related but different from the previously studied Riemannian geometry based on the notion of Fisher-Rao information matrix taken as an expression of an inner product between tangent vectors. This new geometry originated, a least in the statistical community, with the idea of defining the geometry of curved exponential models as embedded in a larger exponential family \cite{cenkov:1982,efron:1975,efron:1978,amari:82}.

Such a theory has been known for a long time in statistical mechanics. The main difference is that R. Fisher and other statisticians of the same period used to think about parsimoniously parameterized models. In contrast, physicists such as Boltzmann and Gibbs used to think in terms of simple relations between statistical observables. The exponential family appears as a model with peculiar invariance properties in statistical mechanics.\footnote{The difference in terminology between Mathematical Statistics and Physics is sometimes confusing. The tutorial \cite{susskind-Hrabovsky:2013} and the textbook \cite{arnold:1989} should be helpful.}

Consider, for \emph{example}, the Hamiltonian $H(y,x) = \frac {y^2}{2m} + V(x)$, $y,x \in \reals$. Let the associated flow be $(t,y,x) \mapsto T_t(y,x) \in \reals^2$, that is,
\begin{align*}
    &\derivby t T_t(y,x) = \omega H(T_t(y,x))\ , \quad \omega H = \left(-\pderivby x H,\pderivby y H\right) \ , \\
    &T_0(y,x) = (y,x) \ . 
\end{align*}
The evolution of a probability measure $f \circ m$, $m(dy,dx) = dydx$, $f > 0$, under the action of the flow is
\begin{equation*}
\mu(t) = (T_t)_\# f \cdot m = f \circ S_t \det J S_t \cdot m =   f_t \cdot m \ , \quad S_t = T_t^{-1} \ ,
\end{equation*}
provided the flow is a global diffeomorphism.  The curve $t \mapsto f_t$ is controlled by the continuity equation,
\begin{equation*}
    \pderivby t f_t + \omega H \cdot \nabla f_t  = 0 \ .
  \end{equation*}
The score is $d\dot\mu_t/d\mu_t = d\log f_t/f_t$ and the continuity equation can be written
\begin{equation*}
    \pderivby t \log f_t + \omega H \cdot \nabla \log f_t  = 0 \ .
  \end{equation*}

In particular, $f_t = f$ if $f$ is a function of $H$. Among all invariant probability densities, the curve $\theta \mapsto p_t = \euler^{\theta H}/Z(\theta)$ represents an evolution in the class of invariant probabilities. The score of the model is
\begin{equation*}
    \derivby t \log \left(\euler^{\theta H}/Z(\theta)\right) = H - \derivby t \log Z(t) = H - \int H \, p_t\,  dm \ .
\end{equation*}
In statistical mechanics, the score is interpreted as the fluctuation of the Hamiltonian.

The differential geometry of the dual collections of H.~Nagaoka \cite[8.4]{amari-nagaoka:2000} naturally follows from the various elements discussed above. In particular, the differential notion of connection can be derived from the notion of parallel transport in an affine setting. In the non-parametric setting, it is convenient to base the affine structure on a variation of the original notion of affine space of \cite{weyl:1952}. Below is a summary of the presentation in \cite{chirco-pistone:2022}.

The word ``affine'' above refers to the geometrical construction of vectors associated with displacement according to classical H.~Weyl's \emph{axioms} of an affine space. Let be given a set $M$ and a real finite-dimensional vector space $V$. A displacement mapping is a mapping
\begin{equation*}
M \times M \ni (P,Q) \mapsto \overrightarrow{PQ} \in V \ ,
\end{equation*}
such that
\begin{enumerate}
\item for each fixed $P$ the partial mapping $s_P \colon Q \mapsto \overrightarrow{PQ}$ is \emph{1-to-1 and onto}, and
\item the parallelogram law, $\overrightarrow{PQ} + \overrightarrow{QR} = \overrightarrow{PR}$, holds true.
\end{enumerate}
The structure $(M,V,\overrightarrow{\phantom{pq}})$ is, by \emph{definition}, the affine space. The corresponding affine manifold is derived from the atlas of charts $s_P \colon M \to V$, $P \in M$. Notice that the change of chart is the choice of a new origin. Such a structure supports a full geometrical development, see \cite{nomizu-sasaki:94}.

Weyl's axioms suggest the following \emph{definition}.

\emph{Let $M$ be a set and  let $B_\mu$, $\mu \in M$, be a family of real topological vector spaces. Let $(\transport {\nu}{\mu})$, $\nu,\mu \in M$ be a family of isomorphism $\transport {\nu} {\mu} \colon B_{\nu} \to B_{\mu}$ satisfying the cocycle condition, 
\begin{description} 
\item[AF0]
 $\transport \mu \mu = I$ and $\transport \nu \rho \transport \mu \nu  = \transport \mu \rho$, where 
 $\transport \nu \mu$ is the transport from $B_\nu$ onto $B_\mu$.
\end{description}
Consider a displacement mapping
\begin{equation*}
\displacement \colon  (\nu,\mu) \mapsto s_\nu(\mu)\in B_{\nu}
\end{equation*}
defined on a subset of the product space $\domof \displacement \subset M \times M$. Assume
\begin{description}\label{desc:affineaxiom}
\item[AF1] \label{item:affineaxiom1}
  For each fixed $\nu$ the partial mapping $M_\nu \ni \mu \mapsto s_\nu(\mu) = \displacement(\nu,\mu)$ is injective. 
\item[AF2] \label{item:affineaxiom2}
  $\displacement (\mu_1,\mu_2) + \transport {\mu_2}{\mu_1} \displacement (\mu_2,\mu_3) = \displacement (\mu_1,\mu_3)$.
  \end{description}
The structure $(M,(B_\mu)_{\mu \in M},(\transport {\nu}{\mu})_{\mu,\nu\in M},\displacement)$ is an affine bundle.}

The affine bundle provides a family of candidates to charts $s_\nu \colon M_\nu \to B_\nu$, $\nu \in M$, to for an atlas. 
Let $\left(M,(B_\mu)_{\mu \in M},(\transport {\nu}{\mu})_{\mu,\nu\in M},\displacement \right)$ be an affine space and assume
\emph{AF3} For each $\nu$, the image set $s_\nu(M_\nu)$ is a neighborhood of 0 in $B_\mu$. That is, its interior $s_\nu(M_\nu)^\circ$ is an open set containing $s_\nu(\nu)=0$. Define the coordinates domains as $U_\nu = s_\nu^{-1} \left(s_\nu(M)^\circ\right)$, so that $(s_\nu, U_\nu, B_\nu)$ is a chart on $M$. Such a chart is said to have \emph{origin} $\nu$. Such charts are compatible, and the resulting manifold
\begin{equation*}
 \M = \left(M,(B_\mu)_{\mu \in M},(\transport {\nu}{\mu})_{\mu,\nu\in M},(s_\mu)_{\mu\in M} \right)   
\end{equation*}is, by \emph{definition}, the affine manifold associated with the given affine bundle.

 Here is our main instance. Consider the exponential transport of \cref{eq:etransport} and define $s_p(q) = \log \frac q p - \int \log \frac q p \ p \cdot dm$. The parallelogram identity is
  \begin{multline*}
  \left(\log \frac q p - \int \log \frac q p \ p \cdot dm\right) + \\  \left(\log \frac r q - \int \log \frac r q \ p \cdot dm - \int  \left(\log \frac r q - \int \log \frac r q \ p \cdot dm\right)  \ dm \right) = \\ \left(\log \frac r p - \int \log \frac r p \ p \cdot dm\right)    
\end{multline*}
The inverse of the chart is easily seen to be 
\begin{equation*}
  s_p^{-1}(u) = \expof{u - K_p(u)} \cdot p \ , \quad K_p(u) = \log \int \euler^u \ p \ dm \ , \quad u \in B_p \ .
\end{equation*}

The dual instance is associated with the mixture transport of \cref{eq:mtransport} and $s_p(q) = \frac q p - 1$. The parallelogram identity is
\begin{equation*}
    \frac q p - 1 + \frac q p \left(\frac r q -1\right) = \frac r p - 1 \ .
\end{equation*}

 Given an affine manifold $\M$, the affine bundle is again an affine manifold. In
 \begin{equation*}
   S\M = \setof{(\mu,v)}{\mu \in M, v \in B_\mu}
 \end{equation*}
 the equation
 \begin{equation*}
   S\M \times S\M \ni ((\nu,u),(\mu,v)) \mapsto (s_\nu(\mu),\transport \mu \nu v) \in B_\nu \times B_\nu
 \end{equation*}
 defines a displacement on the bundle. For each $\nu$ define the chart
 \begin{equation*} 
     s_\nu \colon SM \ni (\mu,v) \mapsto (s_\nu(\mu), \transport \mu \nu v) \in B_\nu \times B_\nu 
 \end{equation*}
 to define the affine bundle $S\M$ as a manifold. Equivalently, we can say that $S\M$ is the bundle with trivialization
 \begin{equation*}
     s_\nu \colon (\mu,v) \mapsto (s_\nu(\mu), \transport \mu \nu v) \ .
 \end{equation*}

The affine bundle is a convenient expression of the tangent bundle of the affine manifold if we \emph{define} the velocity as follows. The velocity of the smooth curve $t \mapsto \gaussdensity(t)$ of the affine manifold $\M$ is the curve $t \mapsto (\gaussdensity(t), \velocity \gaussdensity(t))$ of the affine bundle whose second component is
\begin{equation*}
\velocity \gaussdensity(t) = \lim _{h \to 0} h^{-1} (s_{\gaussdensity(t)}(\gaussdensity(t+h)) = \left. \derivby h s_{\gaussdensity(t)}(\gaussdensity(t+h)) \right\vert_{h=0} \ .
\end{equation*}

By assumption \emph{AF2} applied to the points, the expression in the chart centered at $\nu$ of $\velocity \gaussdensity(t)$ is $\transport {\gaussdensity(t)} \nu \velocity \gaussdensity(t) =  s_\nu(\gaussdensity(t))$. 

For \emph{example}, in the exponential manifold, it holds
\begin{equation*}
    \velocity \gaussdensity(t) = \etransport p {\gaussdensity(t)} \derivby t \left(\log \frac {\gaussdensity(t)} p - \expectat p {\log \frac {\gaussdensity(t)} p}\right) = \derivby t \log \gaussdensity(t) \ ,
\end{equation*}
so that the (affine) velocity in the exponential manifold equals Fisher's score.
  
Let $F$ be a section of the affine manifold, that is, $\mu \mapsto (\mu,F(\mu)) \in \mathcal S \M$. An \emph{integral curve} of the section $F$ is a curve $t \mapsto \gaussdensity(t)$ such that $\velocity \gaussdensity(t) = F(\gaussdensity(t))$. A \emph{flow} of the section $F$ is a mapping
\begin{equation*}
    M \times I \ni (\nu,t) \mapsto \gaussdensity_t(\nu)
\end{equation*}
such that for each $\nu$ the curve $t \mapsto \gaussdensity_t(\nu)$ is an integral curve and $\gaussdensity(0,\nu) = \nu$. 

The following \emph{proposition} gives a characterization of affine geodesics.
  \emph{The following statements are equivalent. 
  \begin{enumerate}
  \item The curve $I \colon t \mapsto \gaussdensity(t)$ is auto parallel, that is, 
  $\velocity \gaussdensity(t) = \transport {\gaussdensity(s)} {\gaussdensity(t)} \velocity \gaussdensity(s)$, $s,t \in I$.
  \item The expression of the curve in each chart is affine.
  \item For all $s,t$  
\begin{equation*} 
  \gaussdensity(t) = S^{-1}_{\gaussdensity(s)}\left((t-s) \velocity \gaussdensity(s)\right) \end{equation*}
\end{enumerate}}


 The acceleration is \emph{defined} as a velocity in the affine bundle. Consider the curve $t \mapsto \mu(t)$ with velocity $t \mapsto \velocity \mu(t)$. The acceleration $t \mapsto \acceleration \mu(t)$ is the velocity $t \mapsto (\mu(t),\velocity \mu(t))$.
\begin{equation*}
     (\velocity \mu(t),\acceleration \mu(t)) = \lim_{h \to 0} h^{-1} s_{\mu(t),\velocity \mu(t)}(\mu(t+h),\velocity \mu(t+h)) \ .
   \end{equation*}
Especially, for all $\mu \in M$,
\begin{equation*}
  \acceleration \mu(t) = \transport \mu {\mu(t)} \derivby t \transport {\mu(t)} \mu \velocity \mu(t) \ .
\end{equation*}
This equation shows that a curve with 0 acceleration is auto-parallel.

In the exponential \emph{example}, the acceleration is computed as follows:
In the exponential case,
\begin{multline*}
  \acceleration p(t) = \etransport {p}{p(t)} \derivby t \etransport {p(t)} p \velocity p(t) = 
  \etransport p {p(t)} \derivby t \left(\frac{\dot p(t)}{p(t)} - \int \frac{\dot p(t)}{p(t)} \ p \ dm\right) = \\
  \frac {\ddot p(t)}{p(t)} - \left(\frac {\dot p(t)}{p(t)}\right)^2 - \int \left(\frac {\ddot p(t)}{p(t)} - \left(\frac {\dot p(t)}{p(t)}\right)^2\right) \ p(t) \ dm = \\
  \frac {\ddot p(t)}{p(t)} - \left(\frac {\dot p(t)}{p(t)}\right)^2 + \int \left(\frac {\dot p(t)}{p(t)}\right)^2 \ p(t) \ dm\ . 
\end{multline*}

\section{Gaussian Orlicz-Sobolev model space}
\label{sec:stat-grad-model}

Above, we have discussed in general terms how to define an affine Banach manifold. We now proceed to instantiate the general formalism into a specific case of Gaussian space. In doing that, the usual toolbox of IG should be extended with other analytical notions. A general reference is \cite{bogachev:2010}. We now restrict our attention to a  particular instance of model Banach space. Precisely, we are going to use the generalization of Lebesgue spaces called Orlicz spaces. General references are the monographs \cite[Ch.~II]{musielak:1983} and \cite[Ch.~VII]{adams-fournier:2003}. The basic technical tools is are the notion of conjugation between convex functions and the analysis of the Gaussian space. I will use my conference paper \cite{pistone:2021PIGTA}.

\subsection{Orlicz spaces}
Assume $\phi \in C[0,+\infty[$ is null at 0, $\phi(0)=0$, strictly  increasing, and $\lim_{u \to +\infty} \phi(u) = +\infty$. Let $\Phi$ be its primitive function with $\Phi(0) = 0$. We call such a function a Young function. The inverse function $\psi = \phi^{-1}$ has the same properties as $\phi$, so that the primitive $\Psi$ with $\Psi(0)=0$ is again a Young function. The couple $(\Phi, \Psi)$, is a couple of conjugate Young functions. The relation is symmetric and we write both $\Psi=\Phi_*$ and $\Phi = \Psi_*$. The Young inequality holds true,
\begin{equation*} 
  \Phi(x) + \Psi(y) \geq xy \ , \quad x,y \geq 0 \ ,
\end{equation*}
and the Legendre equality holds true,
\begin{equation*} 
  \Phi(x) + \Psi(\phi(x)) = x \phi(x) \ , \quad x \geq 0 \ .
\end{equation*}

Here are my notations for specific cases we are going to use:
\begin{align}
&\Phi(x) = {x^p} / p \ , \quad \Psi(y) =  {y^q} / q \ , \quad p,q > 1 \  , \quad 1/p + 1/q = 1 \ ;   \label{eq:cases-1} \\
  &\exp_2(x) = \euler^x - 1 - x \ , \quad (\exp_2)_*(y) = (1+y)\log(1+y) - y \ ; \label{eq:cases-2} \\
  &\cosh_2(x) = \cosh x - 1 \ , \quad (\cosh_2)_*(y) = \int_0^y\sinh^{-1}(v) \ dv \ ; \label{eq:cases-3} \\
&\gausstwo(x) = \expof{\frac12x^2}-1 \ . \label{eq:cases-4}
\end{align}

Given a Young function $\Phi$, and a probability measure $\mu$, the Orlicz space $\orliczof \Phi \mu$ is the Banach space whose closed unit ball is $$\setof{f \in L^0(\mu)}{\int \Phi(\avalof f) \ d\mu \leq 1} \ .$$ The corresponding Minkowski norm is usually called Luxemburg norm,
\begin{equation*}
\normat {\orliczof \Phi \mu} f = \inf \setof{\alpha > 0}{\int \Phi(\alpha^{-1} \avalof f) \ d\mu   \leq 1} \ .
\end{equation*}
The Young inequality provides a separating duality $\scalarat \mu u v = \int uv \ d\mu$ of $\orliczof {\Phi} \mu$ and $\orliczof {\Phi_*} \mu$ such that $\scalarat \mu u v \leq 2 \normat {\orliczof{\Phi} \mu} u  \normat {\orliczof {\Phi_*} \mu} v$. The dual norm is called the Orlicz norm and is equivalent to the Luxembourg norm.

 Domination relation between Young functions implies continuous injection properties for the corresponding Orlicz spaces. We say that $\Phi_2$ \emph{eventually dominates} $\Phi_1$, written $\Phi_1 \prec \Phi_2$, if there is a constant $\kappa$ such that $\Phi_1(x) \leq \Phi_2(\kappa x)$ for all $x$ larger than some $\bar x$. As, in our case, $\mu$ is a probability measure, the continuous embedding $\orliczof {\Phi_2} \mu \to \orliczof {\Phi_1} \mu$ holds if, and only if, $\Phi_1 \prec \Phi_2$. See proof in \cite[Th. 8.2]{adams-fournier:2003}. If $\Phi_1 \prec \Phi_2$, then $(\Phi_2)_* \prec (\Phi_1)_*$. Looking at the examples above, $\exp_2$ \eqref{eq:cases-2} and $\cosh_2$ \eqref{eq:cases-3} are equivalent, they both are eventually dominated by $\operatorname{gauss}_2$ \eqref{eq:cases-4} and eventually dominate all powers \eqref{eq:cases-1}.

A special case occurs when there exists a function $C$ such that $\Phi(ax) \leq C(a) \Phi(x)$ for all $a \geq 0$. This is true, for example, for a power function and in the case of the functions $(\exptwo)_*$ and $(\coshtwo)_*$. In such a case, the conjugate space is the dual space and and bounded functions form a dense set.

The spaces corresponding to case \eqref{eq:cases-1} are ordinary Lebesgue spaces. The cases \eqref{eq:cases-2} and \eqref{eq:cases-3} provide isomorphic Banach spaces, which are of special interest to us as they provide the model spaces for our non-parametric version of IG. In fact, a random variable $u$ belongs to $\orliczof {\exp_2} \mu$ if, and only if, the exponential family $p_\theta \propto \euler^{\theta u}$ is defined in a neighborhood of $\theta = 0$. In the conjugate space, a strictly positive probability density $f$ has finite entropy if, and only if, the random variable $v = f -1$ belongs to $\orliczof {(\exp_2)_*} \mu$.

\subsection{Sub-exponential and sub-Gaussian random variables}
\label{sec:sub-exponential} 
There is another important feature of the class $\orliczof {\cosh_2} \mu$. Such a class coincides with the class of \emph{sub-exponential} random variables, that is, those for which there exist constants $C_1,C_2 > 0$ such that the large deviations admit an exponential bound
\begin{equation*}
  \probat \mu {\avalof{f} \geq t} \leq C_1 \expof{-C_2 t} \ , \quad t \geq 0 \ .
\end{equation*}
Sub-exponential random variables are of special interest in statistical applications because they admit explicit exponential bounds in the Law of Large Numbers. There is a large amount of literature on this subject; see, for example, \cite{buldygin-kozachenko:2000,vershynin:2018-HDP,wainwright:2019-HDS}.

Random variables whose square is sub-exponential are called sub-Gaussian.  For each Young function $\Phi$, the function $\overline \Phi(x) = \Phi(x^2)$ is again a Young function, and $\normat {\orliczof {\overline \Phi} \mu} f \leq \lambda$ if, and only if, $\normat {\orliczof \Phi \mu} {\avalof f ^2} \leq \lambda^2$. We denote the resulting space by $\orliczpof 2 \Phi \mu$. For example, $\gausstwo$ and $\overline{\cosh_2}$ are $\prec$-equivalent , hence the isomorphisn $\orliczof {\gausstwo} \mu \leftrightarrow \orliczpof 2 {\coshtwo} \mu$.

As an application, consider that for each increasing convex $\Phi$ it holds $\Phi(fg) \leq \Phi((f^2+g^2)/2) \leq (\Phi(f^2) + \Phi(g^2))/2$. It follows that when the $\orliczpof \Phi 2 \mu$-norm of $f$ and of $g$ is bounded by one, the $\orliczof \Phi \mu$-norm of $f$, $g$, and $fg$, are all bounded by one. The need to control the product of two random variables in $\orliczof {\coshtwo} \mu$ appears, for example, in the study of the covariant derivatives of the statistical bundle, see \cite{gibilisco-pistone:98,lott:2008calculations,pistone:2018lagrange,chirco-malago-pistone:2022}.

\subsection{Calculus of the Gaussian space}
\label{sec:calc-gauss-space}

From now on, the base probability space is the Gaussian probability space $(\reals^n,\gaussdensity)$, $\gaussdensity(z) = (2\pi)^{n/2} \expof{- \avalof z ^2 / 2}$. Let us recall a few simple facts about the analysis of the Gaussian space, see \cite[Ch.~V]{malliavin:1995}.

Let us denote by $C^k_\text{poly}(\reals^n)$, $k = 0,1,\dots$, the vector space of functions which are differentiable up to order $k$ and which are bounded, together with all derivatives, by a polynomial. This class of functions is dense in $L^2(\gaussdensity)$. For each couple $f,g \in \cpoly 1 n$, we have
\begin{equation*}
  \gaussint {f(x) \ \partial_i g(x)} x =   \gaussint {\delta_i f(x) \ g(x)} x \ ,
\end{equation*}
where the divergence operator $\delta_i$ is defined by $\delta_i f(x) = x_i f(x) - \partial_i f(x)$. Vector case is
\begin{equation*}
  \gaussint {\nabla f(x) \cdot \nabla g(x)} x = \gaussint {f(x) \ \delta \cdot \nabla g(x)} x \ , \quad f,g \in \cpoly 2 n \ ,
\end{equation*}
with $\delta \cdot \nabla g(x) = x \cdot \nabla g(x) - \Delta g(x)$.

Hermite polynomials $H_\alpha = \delta^\alpha 1$ provide an orthogonal basis for $L^2(\gaussdensity)$ such that $\partial_i H_\alpha = \alpha_i H_{\alpha - e_i}$, $e_1$ the $i$-th element of the standard basis of $\reals^n$. In turn, this provides a way to prove that there is a closure of both operator $\partial_i$ and $\delta_i$ on a domain which is an Hilbert subspace of $L^2(\gaussdensity)$. Such a space is denoted by $D^2$ in \cite{malliavin:1995}. Moreover, the closure of $\partial_i$ is the infinitesimal generator of the translation operator, \cite{malliavin:1997,bogachev:2010}. The space $D^2$ is a Sobolev Space with Gaussian weight based on the $L^2$ norm, \cite{adams-fournier:2003}.

\subsection{Exponential statistical bundle}
\label{sec:exponential-bundle}

We refer to \cref{sec:bundle} and \cite{pistone:2013GSI,pistone:2018-IGAIA-IV} for the definition of maximal exponential manifold $\maxexpat \gaussdensity$, and of statistical bundle $S \maxexpat \gaussdensity$. Below we report the results that are necessary in the context of the present paper.

A key result is the proof of the following \emph{proposition}, see \cite{pistone-sempi:95,cena:2002,cena-pistone:2007} and \cite[Th.~4.7]{santacroce-siri-trivellato:2016}.

\emph{For all $p, q \in \maxexpat \gaussdensity$ it holds $q = \euler^{u - K_p(u)} \cdot p$, where $u \in \orliczof {\coshtwo} \gaussdensity$, $\expectat p u = 0$, and $u$ belongs to the interior of the proper domain of the convex function $K_p$. This property is equivalent to any of the following:
\begin{enumerate}\label{portmanteaux}
\item\label{portmanteaux-1} $p$ and $q$ are connected by an open exponential arc;
\item\label{portmanteaux-2} $\Lexp p = \Lexp q$ and the norms are equivalent;
\item\label{portmanteaux-3}
  $p/q \in \cup_{a > 1} L^a(q)$ and $q/p \in \cup_{a>1} L^a(p)$.
\end{enumerate}
}

\Cref{portmanteaux-2} ensures that all the fibers of the statistical bundle, namely $\expfiberat p \gaussdensity$, $p \in \maxexpat \gaussdensity$,  are isomorphic. \Cref{portmanteaux-3} gives a explicit description of the exponential manifold. For \emph{example}, let $p$ be a positive probability density with respect to $\gaussdensity$, and take $q=1$ and $a = 2$. Then a sufficient condition for $p \in \maxexpat \gaussdensity$ is
\begin{equation*}
  \gaussint {p(x)^2} x < \infty \quad \text{and} \quad \gaussint {\frac1{p(x)}} x < \infty \ .
\end{equation*}

\subsection{Orlicz spaces with derivatives}
\label{sec:orlicz-sobolev}
By replacing the $L^2$-norm with a $\coshtwo$-Orlicz norm, a set-up for IG obtains \cite{lods-pistone:2015,pistone:2018-IGAIA-IV}. Precisely, we have exponential families with weakly differentiable densities and a Gaussian reference probability measure. The main outcome is the possibility to discuss topics related to the analytical picture of IG. 

Let us start with a class of inequalities related to the classical Gauss-Poincar\'e inequality,
\begin{equation*} 
\gaussint {\left(f(x) - \gaussint {f(y)} y \right)^2} x \le \gaussint {\avalof{\nabla f(x)}^2} x \ ,
\end{equation*}
where $f \in C^1_\text{poly}(\reals^n)$. See a proof in \cite[1.4]{nourdin-peccati:2012}. In terms of norms, the inequality above is equivalent to $\normat {L^2(\gaussdensity)}{f - \overline f} \leq \normat {L^2(\gaussdensity)} {\avalof{\nabla f}}$, where $\overline f = \gaussint {f(y)} y$.

For \emph{example}, if $p \in C^2_{\text{poly}}$ is a probability density with respect to $\gaussdensity$, then the $\chi^2$-divergence of $P = p \cdot \gaussdensity$ from $\gaussdensity$ is bounded by
\begin{equation*}
 D_{\chi^2}(P \vert \gaussdensity) =  \gaussint {(p(x) - 1)^2} x \leq \gaussint {(\delta\cdot \nabla p(x))^2} x \ .
\end{equation*}

Generalizations of the Gauss-Poincar\'e inequality follow from the properties of the Ornstein-Uhlenbeck (OU) semigroup
\begin{equation*} 
  P_tf(x) = \gaussint {f(\euler^{-t}x + \sqrt{1-\euler^{-2t}}y)} y, \quad t \ge 0,\quad   f \in C^k_\text{poly}(\reals^n) \ ,
\end{equation*}
see \cite[V-1.5]{malliavin:1995} and \cite[1.3]{nourdin-peccati:2012}. Notice that the OU semigroup interpolates between $P_0f = f$ and $P_\infty f = \overline f$. If  $X$, $Y$ are independent standard Gaussian random variables in $\reals^n$, then
\begin{equation*} 
  X_t = \euler^{-t}X+\sqrt{1-\euler^{-2t}}Y, \quad Y_t=\sqrt{1-\euler^{-2t}}X-\euler^{-t}Y
\end{equation*}
are independent standard Gaussian random variables for all $t \ge 0$.
By the change of variable $(X,Y) \to (X_t,Y_t)$ and Jensen's inequality, it follows for each convex $\Phi$ that
\begin{equation*} 
  \gaussint {\Phi(P_tf(x))} x \leq \gaussint {\Phi(f(x))} x \ .
\end{equation*}
That is, for all $t \ge 0$, the mapping $f \mapsto P_tf$ is non-expansive  for the norm of each Orlicz space $\orliczof \Phi \gaussdensity$.

For all $\Phi\colon \reals$ convex and all $f \in C^1_{\text{\normalfont poly}}(\reals^n)$, it holds
\begin{multline} \label{eq:first-version}
 \gaussint {\Phi\left(f(x) - \gaussint {f(y)} y\right)} x  \leq \\ 
 \iint \Phi\left(\frac\pi2 \nabla f(x) \cdot y\right) \ \gaussdensity(x) \gaussdensity(y) \ dx dy = \\
\frac1{\sqrt{2\pi}} \iint \Phi\left(\frac\pi2\avalof{\nabla f(x)}z\right) \ \euler^{-z^2/2} \gaussdensity(x) \ dz dx = \\ \gaussint {\widetilde \Phi\left(\avalof{\nabla f(x)}\right)} x \ , 
\end{multline}
where $\widetilde \Phi$ is the convex function
\begin{equation*}
  \label{eq:tildePhi}
  \widetilde \Phi(a) = \gaussint {\Phi\left(\frac\pi2 az\right)} z \ .
\end{equation*}

The first example of convex function is $\Phi(s) = \euler^{s}$, with $\widetilde \Phi(a) = \expof{\frac {\pi^2 a^2}8}$, so that the inequality applied to $\frac{2\kappa}\pi f$ becomes
\begin{equation} \label{eq:MGF-1}
 \gaussint {\expof{\frac{2\kappa}\pi \left(f(x) - \overline f\right)}} x  \leq \gaussint {\expof{\frac{\kappa^2}2 \avalof{\nabla f(x)}^2}} x
\end{equation}
If the function $f$ is Lipschitz with norm $\kappa^{-1}$, then the RHS is finite.

The first case of bound for Orlicz norms we is the $\Phi(s) = s^{2p}$, $p > 1/2$. In such a case, 
\begin{equation*}
 \widetilde \Phi(a) =  \left(\frac \pi 2\right)^{2p} m(2p)\  a^{2p} \ ,
\end{equation*}
where $m(2p)$ is the $2p$-moment of the standard Gaussian distribution. It follows that
 \begin{equation*}
\normat {L^{2p}(\gaussdensity)} {f - \gaussint {f(y)} y} \leq \frac \pi 2 (m(2p))^{1/2p} \normat {L^{2p}(\gaussdensity)} {\avalof{\nabla f}} \ .
\end{equation*}

The cases $\Phi(a) = a^{2p}$ are special in that we can use them in the proof the multiplicative property $\Phi(ab) = \Phi(a)\Phi(b)$. The argument generalizes to the case where the convex function $\Phi$ is a Young function whose increase is controlled through a function $C$, $\Phi(uv) \leq C(u) \Phi(v)$, and such that there exists a $\kappa > 0$ for which
\begin{equation*}
  \gaussint {C\left(\frac \pi 2 \kappa u\right)} u \leq 1 \ , 
\end{equation*}
so that \cref{eq:tildePhi} becomes
\begin{equation*}
  \widetilde \Phi(\kappa a) = \gaussint {\Phi\left(\frac \pi2 \kappa a z\right)} z \leq \gaussint {C\left(\frac \pi 2 \kappa z\right)} z \ \Phi(a) \leq \Phi(a)\ .
\end{equation*}

By using this bound in \cref{eq:first-version}, we get
\begin{equation*}
   \gaussint {\Phi\left(\kappa \left(f(x) - \gaussint {f(y)} y\right)\right)} x  \leq \gaussint {\Phi\left(\avalof{\nabla f(x)}\right)} x \ .
 \end{equation*}

 Assume now that $\normat {\orliczof {\Phi} {\gaussdensity}} {\avalof{\nabla f}}\leq 1$ so that the LHS does not exceed 1. Then $\kappa \normat {\orliczof {\Phi} {\gaussdensity}} {f - \overline f} \leq 1$, which, in turn, implies the inequality
   \begin{equation*}
     \normat {\orliczof {\Phi} {\gaussdensity}} {f-\overline f} \leq \kappa^{-1} \normat {\orliczof {\Phi} {\gaussdensity}} {\avalof{\nabla f}} \ . 
   \end{equation*}

   It is of particular interest the case of the Young function $\Phi = \cosh-1$, for which there is no such bound. Instead, we use \cref{eq:MGF-1} with $\kappa$ and $-\kappa$ to get
\begin{multline}\label{eq:cosh}
 \gaussint {(\cosh-1)\left(\frac{2\kappa}\pi \left(f(x) - \overline f\right)\right)} x  \leq \\ \gaussint {\operatorname{gauss}_2\left(\kappa \avalof{\nabla f(x)}\right)} x \ . 
\end{multline}
Now, if $\kappa = \normat {\orliczof {\operatorname{gauss}_2} {\gaussdensity}}{\avalof{\nabla f}}^{-1}$, then the LHS is smaller or equal then 1, and hence $2\kappa/\pi \normat {\orliczof {\cosh-1}{\gaussdensity}} {f - \overline f} \leq 1$. It follows that
\begin{equation*}
  \normat {\orliczof {\cosh-1}{\gaussdensity}} {f - \overline f} \leq \frac\pi2 \normat {\orliczof {\operatorname{gauss}_2} {\gaussdensity}}{\avalof{\nabla f}} \ .
\end{equation*}

In the following \emph{proposition}, we summarize the inequalities proved so far.

\emph{There exists constants $C_1$, $C_2(p)$, $C_3$ such that for all $f \in C^1_\text{\emph{poly}}(\reals^n)$ the following inequalities hold:
 \begin{equation}\label{eq:LlogL-bound}
   \normat {L_{(\exp_2)_*}(\gaussdensity)} {f - \gaussint {f(y)} y} \leq C_1 \normat {L_{(\exp_2)_*}(\gaussdensity)}  {\avalof{\nabla f}} \ .
 \end{equation}
 \begin{equation}\label{eq:2k-bound}
   \normat {L^{2p}(\gaussdensity)} {f - \gaussint {f(y)} y} \leq C_2(p) \normat {L^{2p}(\gaussdensity)} {\avalof{\nabla f}} \ , \quad p > 1/2 \ .
 \end{equation}
 \begin{equation}\label{eq:exp-bound}
           \normat {\Lexp \gaussdensity} {f - \gaussint {f(y)} y} \leq C_3 \normat {\orliczof {\operatorname{gauss}_2} \gaussdensity} {\avalof{\nabla f}} \ .
 \end{equation}}

Other equivalent norms could be used in the inequalities above. For
example,
$\orliczof {(\exp_2)_*} \gaussdensity \leftrightarrow \orliczof
{(\cosh-1)_*} \gaussdensity$ and
$\orliczof {\operatorname{gauss}_2} \gaussdensity \leftrightarrow
\orliczpof {\cosh-1} 2 \gaussdensity$.
\label{sec:gener-ornst-uhlenb}

We now consider a further set of inequalities based on the use of infinitesimal generator $- \delta \cdot \nabla$ of the OU semigroup \cite[1.3.7]{nourdin-peccati:2012}.

We have, for all $f \in \cpoly 2 n$, that
\begin{equation}\label{eq:usegenerator-1}
  f(x) - \overline f = - \int_0^{\infty} \derivby t P_t f(x) \ dt = \int_0^{\infty} \delta \cdot \nabla P_t f(x) \ dt \ .
\end{equation}
Note that
\begin{multline*} 
  \nabla P_t f(x) = \nabla \gaussint {f(\euler^{-t}x + \sqrt{1-\euler^{-2t}}y)} y = \\ \euler^{-t} \gaussint {\nabla f(\euler^{-t}x + \sqrt{1-\euler^{-2t}}y)} y = \euler^{-t} P_t \nabla f(x) \ ,
\end{multline*}
so that
\begin{equation*} 
 P_t \delta \cdot \nabla f(x) = \delta \cdot \nabla P_tf(x) = \euler^{-t} \delta \cdot P_t \nabla f(x) \ .
\end{equation*}
Now, \cref{eq:usegenerator-1} becomes
\begin{equation}
  \label{eq:usegenerator-3}
  f(x) - \overline f = \int_0^{\infty} \euler^{-t} \delta \cdot P_t \nabla f(x) \ dt \ .
\end{equation}

As
\begin{equation*}
  \gaussint {\delta \cdot \nabla f(x)} x = 0 \ ,
\end{equation*}
the covariance of $f,g \in \cpoly 0 n$ is
  \begin{multline*}
    \covat \gaussdensity f g = \\ \gaussint{(f(x) - \overline f)g(x)} x = \gaussint{(f(x) - \overline f)(g(x) - \overline g)} x \ .
  \end{multline*}

It follows that for all $f, g \in \cpoly 2 n$ we derive from \cref{eq:usegenerator-3} 
\begin{equation}\label{eq:OU-selfadjoint}
  \covat \gaussdensity f g = \int_0^{\infty} \euler^{-t} \gaussint {P_t \nabla f(x) \cdot \nabla g(x)} x \ dt \ . 
 \end{equation}

We use here a result of \cite[Prop. 5]{pistone:2018-IGAIA-IV}. Let $\avalof \cdot _1$ and $\avalof \cdot _2$ be two norms on $\reals^n$, such that $\avalof{x \cdot y} \le \avalof x _1 \avalof y _2$. For a Young function $\Phi$, consider the norm of $\orliczof {\Phi} \gaussdensity$ and the conjugate space endowed with the dual norm, 
 \begin{equation*}
   \normat {\orliczof {\Psi,*} \gaussdensity} f = \sup\setof{\int fg\  \gaussdensity}{\int \Phi(g)\ \gaussdensity \leq 1} \ .
 \end{equation*}

The following \emph{proposition} includes the standard Poincar\'e case provided $\Phi(u) = u^2/2$. 

\emph{Given a couple of conjugate Young function $\Phi$, $\Psi$, and norms $\avalof\cdot_1$, $\avalof\cdot _2$ on $\reals^n$ such that $x \cdot y \leq \avalof x _1 \avalof y _2$, $x,y \in \reals^n$, for all $f,g \in C^1_\text{\emph{poly}}(\reals^n)$, it holds
\begin{equation*}
  \avalof{\covat \gaussdensity  f g} \le \normat {\orliczof {\Phi} {\gaussdensity}}{\avalof{\nabla f}_1} \normat {\orliczof {\Psi,*} {\gaussdensity}} {\avalof{\nabla g}_2}  \ .  
\end{equation*}}

The case of our interest here is $\Phi = \cosh-1$, $\Psi=(\cosh-1)_*$. As $(\cos-1)_* \prec (\cosh-1)$, it follows, in particular, that $\covat \gaussdensity f f$ is bounded by constant times $\normat {\orliczof {\cosh-1} \gaussdensity} {\avalof{\nabla f}}^2$.

\subsection{Orlicz-Sobolev space with Gaussian weight}
\label{sec:gauss-orlicz-sobolev}
A reasonable option for our model space is to assume densities $f = \euler^{u - K_\gaussdensity(u)} \cdot \gaussdensity$ in the Gaussian maximal exponential family, $f \in \maxexpat \gaussdensity$, and, moreover, assume differentiability in the form $u \in \orliczpof 2 {\coshtwo} \gaussdensity = \orliczof {\gausstwo} \gaussdensity$, that is, $u^2 \in \orliczof {\coshtwo} \gaussdensity$, see \cref{sec:sub-exponential}.

Precisely, the exponential and the mixture Orlicz-Sobolev-Gauss (OSG) spaces of interest are, respectively,
\begin{align}
\Wexptwo \gaussdensity &= \setof{v \in \Lexp \gaussdensity}{\partial_j v \in \orliczof {\gausstwo} \gaussdensity} \label{eq:OrSobExp} \ , \\
\Wlogtwo \gaussdensity &= \setof{\eta \in \LlogL \gaussdensity}{\partial_j f \in \orliczof {\gausstwo_*} \gaussdensity} \label{eq:OrSobLog} \ ,
\end{align}
where $\partial_j$, $ j = 1, \dots, n$, is the partial derivative in the sense of distributions.

 As $\phi \in \Cinfcs{\reals^n}$ implies $\phi \gaussdensity \in \Cinfcs{\reals^n}$, for each $f \in \Wlogtwo \gaussdensity$ it holds
\begin{equation*} \scalarat \gaussdensity {\partial_jf} \phi = \scalarof{\partial_j f}{\phi \gaussdensity} = - \scalarof{f}{\partial_j(\phi \gaussdensity)} = \scalarof {f}{\gaussdensity(u_j - \partial_j)\phi} = \scalarat \gaussdensity f {\delta_j \phi} \ ,
\end{equation*} with $\delta_j \phi = (u_j - \partial_j)\phi$. Here, the \emph{Stein operator} $\delta_i$ acts on $\Cinfcs{\reals^n}$ \cite{brezis:2011fasspde}.

The meaning of both operators $\partial_j$ and $\delta_j = (u_j - \partial_j)$ when acting on square-integrable random variables of the Gaussian space is well known. Still, here we are specifically interested in the action on OSG spaces. Let us denote by $\Cinfp {\reals^n}$ the space of infinitely differentiable functions with polynomial growth of all derivatives. Polynomial growth implies the existence of $\gaussdensity$-moments of all derivatives, hence $\Cinfp{\reals^n} \subset \Wlogtwo \gaussdensity$. If $f \in \Cinfp{\reals^n}$, then the distributional derivative and the ordinary derivative are equal and moreover $\delta_j f \in \Cinfp{\reals^n}$. For each $\phi \in \Cinfcs{\reals^n}$ we have $\scalarat \gaussdensity {\phi}{\delta_j f} = \scalarat \gaussdensity {\partial_j \phi} f$.

The OSG spaces $\Wexptwo \gaussdensity$ and $\Wlogtwo \gaussdensity$ are both Banach spaces \cite[Sec. 10]{musielak:1983}. The norm is the graph norm,
\begin{gather*}
\normat {\Wexptwo \gaussdensity} v = \normat {\Lexp \gaussdensity} v + \sum_{j=1}^n \normat {\orliczof \gausstwo \gaussdensity} {\partial_j v} \ , \\
\normat {\Wlogtwo \gaussdensity} \eta = \normat {\LlogL \gaussdensity} \eta + \sum_{j=1}^n \normat {\orliczof {\gausstwo_*} \gaussdensity} {\partial_j \eta} \ .
\end{gather*}
In the cases of null integral, \cref{eq:exp-bound} shows that the second term only provides an equivalent norm for $\Wexptwo \gaussdensity$.

We review some relations between OSG spaces and Sobolev spaces without weight \cite{adams-fournier:2003} in the following \emph{proposition}. For each ball radius $R > 0$,
\begin{equation*} (2\pi)^{- \frac n2} \ge \gaussdensity(x) \ge \gaussdensity(x) (\avalof x < R) \ge (2\pi)^{- \frac n2} \euler^{-\frac {R^2}2} (\avalof x < R), \quad x \in \reals^n.
\end{equation*}

\emph{Let $\Omega_R$ denote the open sphere of radius $R>0$ and consider the restriction $u \mapsto u_R$ of $u$ to $\Omega_R$.
    \begin{enumerate}
    \item We have the continuous mappings
      \begin{equation*} \Wexp{\reals^n} \subset \Wexp \gaussdensity \rightarrow W^{1,p}(\Omega_R), \quad p \ge 1.
            \end{equation*}
          \item We have the continuous mappings
            \begin{equation*} W^{1,p}(\reals^n) \subset \WlogL{\reals^n} \subset \WlogL \gaussdensity \rightarrow W^{1,1}(\Omega_R), \quad p > 1.
                    \end{equation*}
                    \item \label{item:embeddings3}Each $u \in \Wexp \gaussdensity$ is a.s. H\"older of all orders on each $\overline\Omega_R$ and hence a.s. continuous. The restriction $\Wexp \gaussdensity \to C(\overline\Omega_R)$ is compact.
                    \end{enumerate}}

For \emph{example},if $q = \euler^{v-K_\gaussdensity(v)} \cdot \gaussdensity$ and $p = \euler^{u-K_\gaussdensity(u)} \cdot \gaussdensity$ with $q,p \in \maxexp \gaussdensity$ and $v,u \in \Wexptwo \gaussdensity$, the the Hyv\"arinen divergence is
\begin{equation*}
\label{eq:hyvarinen}
\KH p q = \frac12 \gaussint {\normof{\nabla (u - v)}^2 p(x)} x < + \infty
\end{equation*}
because $\nabla (u-v) \in \Lexp \gaussdensity = \Lexp p \subset L^2(p \cdot \gaussdensity)$.

\section{Conclusion}
\label{sec:discussion}
In this final section, I suggest a few applications of my infinite-dimensional setup of professor Amari's ideas I have considered recently, starting from \cite{lods-pistone:2015} and \cite{pistone:2018-IGAIA-IV}.
I will conclude by mentioning a few other topics, for which it is probably possible to extend rigorous results from the finite state space to the Gaussian space. 

\subsection{Sub-exponential random variables}
\label{sec:rand-vari-with}
Let $f \in \cpoly 2 n$ and assume $f$ is globally Lipschitz, that is, $\avalof{\nabla f(x)} \leq \normat {\operatorname{Lip}(\reals^n)} f$, where $\normat {\operatorname{Lip}(\reals^n)} f$ is the Lipschitz semi-norm. It follows from \cref{eq:MGF-1} that $f \in \orliczof {\coshtwo} \gaussdensity$ and the norm admits a computable bound. If $p$ is any probability density of the maximal exponential model of $\gaussdensity$, that is, it is connected to 1 by an open exponential arc, then the proposition in \cref{sec:exponential-bundle} implies that $f \in \orliczof {\coshtwo} p$, that is, $f$ is sub-exponential under the distribution $P = p \cdot \gaussdensity$. If the sequence $(X_n)_{n=1}^\infty$ is an independent sample of $p \cdot \gaussdensity$, then the sequence of sample means will converge,
\begin{equation*}
  \lim_{n \to \infty} \frac 1n \sum_{j=1}^n f(X_j) = \gaussint {f(x) \ p(x)} x \ ,
\end{equation*}
with an exponential bound on the tail probability. See \cite[2.8]{vershynin:2018-HDP} and \cite{siri-trivellato:2021}.

\subsection{Hyv\"arinen divergence}
\label{sec:hyvarinen-divergence}
I adapt \cite{hivarinen:2005}, \cite{parry-dawid-lauritzen:2012}, and \cite[13.6.2]{amari:2016} to my Gaussian case. Consider the Hyv\"arinen divergence of \cref{eq:hyvarinen} in the Gaussian case, that is, $P = p \cdot \gaussdensity$ and $Q = q \cdot \gaussdensity$. As a function of $q$, the divergence is
\begin{multline*}
  q \mapsto \KH {q \cdot \gaussdensity}{p \cdot \gaussdensity} = \frac12 \gaussint {\avalof{\nabla \log p(x)}^2 p(x)} x + \\ \frac12 \gaussint {\avalof{\nabla \log q(x)}^2 p(x)} x - \gaussint {\nabla \log p(x) \cdot \nabla \log q(x) \ p(x) } x \ , 
\end{multline*}
where the first term does not depend on $q$ and the second term is a $p\cdot\gaussdensity$-expectation. As $\nabla \log p = p^{-1} \nabla p$,  the third term equals
\begin{equation*}
  - \gaussint {\delta \cdot \nabla \log q(x) \ p(x)} x \ ,
\end{equation*}
which is again a $p$-expectation. To minimize the Hyv\"arinen divergence we must minimize the $p$-expected value of the local score
\begin{equation*}
  S(q,x) = \frac12 \avalof{\nabla \log q(x)}^2 - \delta \cdot \nabla \log q(x)  
\end{equation*}

If $p$ and $q$ belong to the maximal exponential model of $\gaussdensity$, then $q = \euler^{u - K(u)}$ with $u \in \orliczof {\coshtwo} \gaussdensity$ and $\gaussint {u(x)} x = 0$. The local score becomes $\frac12 \avalof{\nabla u}^2 - \delta \cdot \nabla u$. To compute the $p$-expected value of the score with an independent sample of $p \cdot \gaussdensity$, we have an interest to assume that the score is in $\orliczof {\coshtwo} \gaussdensity$, because this assumption implies the good convergence of the empirical means for all $p$. Assume, for \emph{example}, $\nabla u \in \orliczpof 2 {\coshtwo} \gaussdensity = \orliczof {\gausstwo} \gaussdensity$. This implies directly $\avalof{\nabla u}^2 \in \orliczof {(\cosh -1)} \gaussdensity$. Moreover, we must assume that the $\orliczof {\coshtwo} \gaussdensity$-norm of $\delta \cdot \nabla u$ is finite. Under such assumptions, one hopes that the minimization of a suitable model of the sample expectation of the Hyv\"arinen score is consistent. 

\subsection{Otto's metric}
\label{sec:ottos-riem-metr}
This metric was originally defined in \cite{otto:2001}. Let be given in the maximal exponential model of $\gaussdensity$, $p \in \maxexpat \gaussdensity$, and let $f$ and $g$ be given in the $p$-fiber of the statistical statistical bundle, that is, $f,g \in \Wexptwo p = \Wexptwo \gaussdensity$ and $\gaussint {f(x)\, p} x = \gaussint {g(x)\, p} x = 0$. Otto's inner product is
\begin{multline*}
(f,g) \mapsto \langle\langle f,g \rangle\rangle _ p= \gaussint {\nabla f(x) \cdot \nabla g(x) \  p(x)} x = \\ \gaussint {f(x) \ \delta \cdot (p(x) \nabla g(x))} x = \scalarat p {f} {\delta \cdot (p \nabla g)} \ . 
\end{multline*}
The LHS is well defined and continuous if $\nabla f, \nabla g \in \orliczpof 2 {\coshtwo}  \gaussdensity$, because, in such a case, $\avalof{\nabla f}^2, \avalof{\nabla g}^2 \in \orliczof {\coshtwo} \gaussdensity = \orliczof {\coshtwo} p$. The RHS, if defined,  is an inner product in $\orliczof {\coshtwo} \gaussdensity$. Note that the mapping $g \mapsto \delta \cdot (p \nabla g)$ is 1-to-1 if $g$ because by $\gaussint {g(x)p(x)} x = 0$. The inverse of this mapping provides the natural gradient of Otto's inner product in the sense of \cite{amari:1998natural} and \cite{li-montufar:2018}.  

\subsection{Boltzmann equation}
The space-homogeneous Boltzmann operator with angular collision kernel $B(z,x) = \avalof{x'z}$ is discussed, for example, in \cite{villani:2002review}. I briefly show below how to use Gaussian Orlicz spaces in this context, see \cite{pistone:2013entropy} and \cite{lods-pistone:2015}.

The mixture bundle $\mixbundleat \gaussdensity$ is the set of all couples $(f,\eta)$ with $f \in \maxexpat \gaussdensity$, $\eta \in \orliczof {\coshtwo_*} \gaussdensity$, and $\gaussint {\eta(x) \, f(x)} x = 0$. The $^*$-notation recalls that the dual of $\orliczof {\coshtwo_*} \gaussdensity$ is $\orliczof {\coshtwo} \gaussdensity$. The Boltzmann operator is the mapping $\maxexpat \gaussdensity \mapsto Q(f)$ with
\begin{multline*}
Q(f)(v) = \\
\int_{\reals^3}\int_{S^2} (f(v - x {x'}(v - w))f(w + x {x'}(v - w))-f(v)f(w))\avalof{{x'}(v - w)}\ \sigma(dx)\ dw \ ,    
\end{multline*}
where $x'$ is the transpose of the column vector $x$, $S_2$ is the unit sphere of $\reals^3$, and $\sigma$ is the uniform probability on $S_2$. One can prove that $f \mapsto Q(f)/f$ is a section of the mixture bundle. The Boltzmann equation can be seen as the equation $\velocity f = Q(f)/f$. 

The smoothness of the Boltzmann section follows from a superposition of operators:
\begin{enumerate}
\item Product: $\maxexp {f_0} \ni f \mapsto f\otimes f \in \maxexp{f_0\otimes f_0}$;
\item Interaction: $\maxexp{f_0\otimes f_0} \ni f\otimes f \mapsto g = B f\otimes f \in \maxexp{f_0\otimes f_0}$;
\item Conditioning: $\maxexp{f_0\otimes f_0} \ni g \mapsto \integrald {S_2}{g\circ A_x}{\sigma(dx)} \in \maxexp{f_0\otimes f_0}$;
\item Marginalization.
\end{enumerate}

There is a weak form of the Boltzmann section. Let $v,w$ be a couple of velocities before collision and let us denote by $(v_x,w_x)$ the velocities after collision, see \cite{villani:2002review}. For $f \in \maxexpat \gaussdensity$ and $g \in \orliczof {\coshtwo} \gaussdensity$, define the operator $A$ with
\begin{equation*}
    Ag(v,w) = \integrald{S^2}{\frac12(g(v_x)+g(w_x))}{\sigma(dx)} - \frac12(g(v)+g(w)) \ .
  \end{equation*}
$Ag$ is in $\orliczof {\coshtwo} {\gaussdensity^{\otimes 2}}$ and $\scalarat f g {Q(f)/f} = \expectat {f\otimes f} {Ag}$.

\subsection{Other possible applications}
The two cases below are open suggestions.
\begin{enumerate}
\item
The transport problem \cite{peyre-cuturi:2019} is discussed from the point of view compatible with IG in \cite{malago-montrucchio-pistone:2018,montrucchio-pistone:2021}. But it does not fit the strictly positive probability densities assumption. However, suppose the given margins belong to the exponential manifold. In that case, it is possible to consider the exponential sub-bundle with the given margins and discuss the gradient flow of the given optimization problem. This is done in a particular finite case in \cite{pistone:2021LNCS}.    
\item The dual couple $\Wexptwo \gaussdensity$ and $\Wlogtwo \gaussdensity$ define a couple of dual bundles on the set $\maxexpat \gaussdensity$ where velocity, acceleration, moment, gradient, Hessian, are all well defined from the affine structure. A real function on the exponential bundle is a Lagrangian function, and a real function on the mixture bundle is a Hamiltonian function. The conjugation relation holds, and the mechanic's equations provide a dynamic picture of the statistical bundle. For example, one can consider the Lagrangian function, where entropy takes the role of the potential energy, and Fisher's metric takes the role of the kinetic energy. This was done in the finite case in \cite{pistone:2018lagrange} and \cite{chirco-malago-pistone:2022}.
\end{enumerate}

\subsection*{Acknowledgments}
The author is partially supported by de Castro Statistics and Collegio Carlo Alberto. He is a member of GNAMPA-INDAM and former faculti of DISMA Politecnico di Torino.
\bibliographystyle{amsplain}

\begin{thebibliography}{10}

\bibitem{adams-fournier:2003}
Robert~A. Adams and John J.~F. Fournier, \emph{Sobolev spaces}, second ed.,
  Pure and Applied Mathematics (Amsterdam), vol. 140, Elsevier/Academic Press,
  Amsterdam, 2003. \MR{2424078 (2009e:46025)}

\bibitem{amari:82}
Shun-Ichi Amari, \emph{{Differential Geometry of Curved Exponential
  Families-Curvatures and Information Loss}}, The Annals of Statistics
  \textbf{10} (1982), no.~2, 357 -- 385.

\bibitem{amari:87dual}
Shu{n-i}chi Amari, \emph{Dual connections on the {H}ilbert bundles of
  statistical models}, Geometrization of statistical theory (Lancaster, 1987)
  (Lancaster) (C.~T.~J. Dodson, ed.), ULDM Publ., 1987, pp.~123--151.

\bibitem{amari:1998natural}
Shun-Ichi Amari, \emph{Natural gradient works efficiently in learning}, Neural
  Computation \textbf{10} (1998), no.~2, 251--276.

\bibitem{amari:2016}
Shun-ichi Amari, \emph{Information geometry and its applications}, Applied
  Mathematical Sciences, vol. 194, Springer, Tokyo], 2016. \MR{3495836}

\bibitem{amari-kumon:1988-AS}
Shun-ichi Amari and Masayuki Kumon, \emph{Estimation in the presence of
  infinitely many nuisance parameters---geometry of estimating functions}, Ann.
  Statist. \textbf{16} (1988), no.~3, 1044--1068.

\bibitem{amari-nagaoka:2000}
Shun-ichi Amari and Hiroshi Nagaoka, \emph{Methods of information geometry},
  Translations of Mathematical Monographs, vol. 191, American Mathematical
  Society, Providence, RI; Oxford University Press, Oxford, 2000, Translated
  from the 1993 Japanese original by Daishi Harada. \MR{1800071}

\bibitem{arnold:1989}
V.~I. Arnold, \emph{Mathematical methods of classical mechanics}, Graduate
  Texts in Mathematics, vol.~60, Springer-Verlag, New York, 1989, Translated
  from the 1974 Russian original by K. Vogtmann and A. Weinstein, Corrected
  reprint of the second (1989) edition. \MR{1345386}

\bibitem{ay-jost-le-schwachhofer:2017ig}
Nihat Ay, J\"urgen Jost, H\^ong~V\^an L\^e, and Lorenz Schwachh\"ofer,
  \emph{Information geometry}, Ergebnisse der Mathematik und ihrer
  Grenzgebiete. 3. Folge., vol.~64, Springer, Cham, 2017. \MR{3701408}

\bibitem{bauer-bruveris-michor:2016}
Martin Bauer, Martins Bruveris, and Peter~W. Michor, \emph{Uniqueness of the
  {F}isher-{R}ao metric on the space of smooth densities}, Bull. Lond. Math.
  Soc. \textbf{48} (2016), no.~3, 499--506. \MR{3509909}

\bibitem{bogachev:2010}
Vladimir~I. Bogachev, \emph{Differentiable measures and the {M}alliavin
  calculus}, Mathematical Surveys and Monographs, vol. 164, American
  Mathematical Society, Providence, RI, 2010. \MR{2663405}

\bibitem{bourbaki:71variete}
Nicolas Bourbaki, \emph{Variétés differentielles et analytiques. fascicule de
  résultats / paragraphes 1 à 7}, Éléments de mathématiques, no. XXXIII,
  Hermann, Paris, 1971.

\bibitem{brezis:2011fasspde}
Haim Brezis, \emph{Functional analysis, {S}obolev spaces and partial
  differential equations}, Universitext, Springer, New York, 2011. \MR{2759829
  (2012a:35002)}

\bibitem{brown:86}
Lawrence~D. Brown, \emph{Fundamentals of statistical exponential families with
  applications in statistical decision theory}, IMS Lecture Notes. Monograph
  Series, no.~9, Institute of Mathematical Statistics, Hayward, 1986.
  \MR{MR882001 (88h:62018)}

\bibitem{buldygin-kozachenko:2000}
V.~V. Buldygin and Yu.~V. Kozachenko, \emph{Metric characterization of random
  variables and random processes}, Translations of Mathematical Monographs,
  vol. 188, American Mathematical Society, Providence, RI, 2000, Translated
  from the 1998 Russian original by V. Zaiats. \MR{1743716}

\bibitem{cena:2002}
Alberto Cena, \emph{Geometric structures on the non-parametric statistical
  manifold}, Ph.D. thesis, Universit\`a degli Studi di Milano, 2002.

\bibitem{cena-pistone:2007}
Alberto Cena and Giovanni Pistone, \emph{Exponential statistical manifold},
  Ann. Inst. Statist. Math. \textbf{59} (2007), no.~1, 27--56. \MR{MR2396032
  (2009b:62011)}

\bibitem{chirco-malago-pistone:2022}
Goffredo Chirco, Luigi Malag{\`{o}}, and Giovanni Pistone, \emph{Lagrangian and
  {H}amiltonian dynamics for probabilities on the statistical bundle},
  International Journal of Geometric Methods in Modern Physics (2022).

\bibitem{chirco-pistone:2022}
Goffredo Chirco and Giovanni Pistone, \emph{Dually affine information geometry
  modeled on a banach space}, 2022.

\bibitem{efron:1975}
Bradley Efron, \emph{Defining the curvature of a statistical problem (with
  applications to second order efficiency)}, Ann. Statist. \textbf{3} (1975),
  no.~6, 1189--1242, With a discussion by C. R. Rao, Don A. Pierce, D. R. Cox,
  D. V. Lindley, Lucien LeCam, J. K. Ghosh, J. Pfanzagl, Niels Keiding, A. P.
  Dawid, Jim Reeds and with a reply by the author. \MR{MR0428531 (55 \#1552)}

\bibitem{efron:1978}
\bysame, \emph{The geometry of exponential families}, Ann. Statist. \textbf{6}
  (1978), no.~2, 362--376. \MR{57 \#10890}

\bibitem{efron-hastie:2016}
Bradley Efron and Trevor Hastie, \emph{Computer age statistical inference},
  Institute of Mathematical Statistics (IMS) Monographs, vol.~5, Cambridge
  University Press, New York, 2016, Algorithms, evidence, and data science.
  \MR{3523956}

\bibitem{gibilisco-pistone:98}
Paolo Gibilisco and Giovanni Pistone, \emph{Connections on non-parametric
  statistical manifolds by {O}rlicz space geometry}, IDAQP \textbf{1} (1998),
  no.~2, 325--347. \MR{1 628 177}

\bibitem{hivarinen:2005}
Aapo Hyv\"{a}rinen, \emph{Estimation of non-normalized statistical models by
  score matching}, J. Mach. Learn. Res. \textbf{6} (2005), 695--709.
  \MR{2249836}

\bibitem{kass-vos:1997}
Robert~E. Kass and Paul~W. Vos, \emph{Geometrical foundations of asymptotic
  inference}, Wiley Series in Probability and Statistics: Probability and
  Statistics, John Wiley \& Sons, New York, 1997.

\bibitem{lang:1995}
Serge Lang, \emph{Differential and {R}iemannian manifolds}, third ed., Graduate
  Texts in Mathematics, vol. 160, Springer-Verlag, New York, 1995.
  \MR{96d:53001}

\bibitem{le:2022arXiv}
H\^ong~V\^an L\^e, \emph{Natural differentiable structures on statistical
  models and the fisher metric}, 2022.

\bibitem{li-montufar:2018}
Wuchen Li and Guido Mont{\'u}far, \emph{Natural gradient via optimal
  transport}, Information Geometry \textbf{1} (2018), no.~2, 181--214.

\bibitem{lods-pistone:2015}
Betrand Lods and Giovanni Pistone, \emph{Information geometry formalism for the
  spatially homogeneous {B}oltzmann equation}, Entropy \textbf{17} (2015),
  no.~6, 4323--4363.

\bibitem{lott:2008calculations}
John Lott, \emph{Some geometric calculations on {W}asserstein space}, Comm.
  Math. Phys. \textbf{277} (2008), no.~2, 423--437. \MR{2358290}

\bibitem{malago-montrucchio-pistone:2018}
Luigi Malag{\`o}, Luigi Montrucchio, and Giovanni Pistone, \emph{Wasserstein
  riemannian geometry of gaussian densities}, Information Geometry \textbf{1}
  (2018), no.~2, 137--179.

\bibitem{malliavin:1995}
Paul Malliavin, \emph{Integration and probability}, Graduate Texts in
  Mathematics, vol. 157, Springer-Verlag, New York, 1995, With the
  collaboration of H\'el\'ene Airault, Leslie Kay and G\'erard Letac, Edited
  and translated from the French by Kay, With a foreword by Mark Pinsky.
  \MR{MR1335234 (97f:28001a)}

\bibitem{malliavin:1997}
\bysame, \emph{Stochastic analysis}, Grundlehren der Mathematischen
  Wissenschaften [Fundamental Principles of Mathematical Sciences], vol. 313,
  Springer-Verlag, Berlin, 1997. \MR{MR1450093 (99b:60073)}

\bibitem{montrucchio-pistone:2021}
Luigi Montrucchio and Giovanni Pistone, \emph{Kantorovich distance on finite
  metric spaces: Arens--eells norm and cut norms}, Information Geometry (2021).

\bibitem{musielak:1983}
Julian Musielak, \emph{Orlicz spaces and modular spaces}, Lecture Notes in
  Mathematics, vol. 1034, Springer-Verlag, Berlin, 1983.

\bibitem{nomizu-sasaki:94}
Katsumi Nomizu and Takeshi Sasaki, \emph{Affine differential geometry: geometry
  of affine immersions}, Cambridge Tracts in Mathematics, no. 111, Cambridge
  University Press, Cambridge, 1994.

\bibitem{nourdin-peccati:2012}
Ivan Nourdin and Giovanni Peccati, \emph{Normal approximations with {M}alliavin
  calculus. from stein's method to universality}, Cambridge Tracts in
  Mathematics, vol. 192, Cambridge University Press, Cambridge, 2012.

\bibitem{otto:2001}
Felix Otto, \emph{The geometry of dissipative evolution equations: the porous
  medium equation}, Comm. Partial Differential Equations \textbf{26} (2001),
  no.~1-2, 101--174. \MR{2002j:35180}

\bibitem{parry-dawid-lauritzen:2012}
Matthew Parry, A.~Philip Dawid, and Steffen Lauritzen, \emph{Proper local
  scoring rules}, Ann. Statist. \textbf{40} (2012), no.~1, 561--592.
  \MR{3014317}

\bibitem{peyre-cuturi:2019}
Gabriel Peyr\'e and Marco Cuturi, \emph{Computational optimal transport},
  Foundations and Trends in Machine Learning \textbf{11} (2019), no.~5--6,
  355--607, arXiv:1803.00567.

\bibitem{pistone:2013entropy}
Giovanni Pistone, \emph{Examples of the application of nonparametric
  information geometry to statistical physics}, Entropy \textbf{15} (2013),
  no.~10, 4042--4065. \MR{3130268}

\bibitem{pistone:2013GSI}
\bysame, \emph{Nonparametric information geometry}, Geometric science of
  information (Frank Nielsen and Fr\'ed\'eric Barbaresco, eds.), Lecture Notes
  in Comput. Sci., vol. 8085, Springer, Heidelberg, 2013, First International
  Conference, GSI 2013 Paris, France, August 28-30, 2013 Proceedings,
  pp.~5--36. \MR{3126029}

\bibitem{pistone:2018-IGAIA-IV}
\bysame, \emph{Information geometry of the {G}aussian space}, Information
  geometry and its applications, Springer Proc. Math. Stat., vol. 252,
  Springer, Cham, 2018, pp.~119--155. \MR{3876116}

\bibitem{pistone:2018lagrange}
\bysame, \emph{Lagrangian function on the finite state space statistical
  bundle}, Entropy \textbf{20} (2018), no.~2, 139.

\bibitem{pistone:2021PIGTA}
\bysame, \emph{Information geometry of smooth densities on the gaussian space:
  Poincar{\'e} inequalities}, pp.~1--17, Springer International Publishing,
  Cham, 2021.

\bibitem{pistone:2021LNCS}
\bysame, \emph{Statistical bundle of the transport model}, Geometric science of
  information, Lecture Notes in Comput. Sci., vol. 12829, Springer, Cham,
  [2021] \copyright 2021, pp.~752--759. \MR{4424383}

\bibitem{pistone-sempi:95}
Giovanni Pistone and Carlo Sempi, \emph{An infinite-dimensional geometric
  structure on the space of all the probability measures equivalent to a given
  one}, Ann. Statist. \textbf{23} (1995), no.~5, 1543--1561. \MR{97j:62006}

\bibitem{santacroce-siri-trivellato:2016}
Marina Santacroce, Paola Siri, and Barbara Trivellato, \emph{New results on
  mixture and exponential models by {O}rlicz spaces}, Bernoulli \textbf{22}
  (2016), no.~3, 1431--1447. \MR{3474821}

\bibitem{siri-trivellato:2021}
Paola Siri and Barbara Trivellato, \emph{Robust concentration inequalities in
  maximal exponential models}, Statistics {\&} Probability Letters \textbf{170}
  (2021), 109001.

\bibitem{susskind-Hrabovsky:2013}
Leonard Susskind and George Hrabovsky, \emph{The theoretical minimum: What you
  need to know to start doing physics}, Basic Books, New York, 2013.

\bibitem{vershynin:2018-HDP}
Roman Vershynin, \emph{High-dimensional probability: an introduction with
  applications in data science}, Cambridge Series in Statistical and
  Probabilistic Mathematics, vol.~47, Cambridge University Press, Cambridge,
  2018, With a foreword by Sara van de Geer. \MR{3837109}

\bibitem{villani:2002review}
C\'dric Villani, \emph{A review of mathematical topics in collisional kinetic
  theory}, Handbook of mathematical fluid dynamics, {V}ol. {I}, North-Holland,
  Amsterdam, 2002, pp.~71--305. \MR{1942465 (2003k:82087)}

\bibitem{wainwright:2019-HDS}
Martin~J. Wainwright, \emph{High-dimensional statistics: a non-asymptotic
  viewpoint}, Cambridge Series in Statistical and Probabilistic Mathematics,
  Cambridge University Press, Cambridge, 2019.

\bibitem{weyl:1952}
Hermann Weyl, \emph{Space {T}ime {M}atter}, Dover, New York, 1952 (eng),
  translation of the 1921 RAUM ZEIT MATERIE.

\bibitem{cenkov:1982}
N.~N. Čencov, \emph{Statistical decision rules and optimal inference},
  Translations of Mathematical Monographs, vol.~53, American Mathematical
  Society, Providence, R.I., 1982, Translation from the Russian edited by Lev
  J. Leifman. \MR{MR645898 (83g:62004)}

\end{thebibliography}

\providecommand{\bysame}{\leavevmode\hbox to3em{\hrulefill}\thinspace}
\providecommand{\MR}{\relax\ifhmode\unskip\space\fi MR }
\providecommand{\MRhref}[2]{%
  \href{http://www.ams.org/mathscinet-getitem?mr=#1}{#2}
}
\providecommand{\href}[2]{#2}

\end{document}